\documentstyle[12pt]{article}

\normalbaselines
\begin{document}
\begin{center}
\vspace*{1.0cm}

{\LARGE{\bf Nonstandard q-deformation of the universal enveloping algebra
${\bf U}({\bf so}_{\bf n})$ }}

\vskip 1.5cm

{\large {\bf A. U. Klimyk }}

\vskip 0.5 cm

Bogolyubov Institute for Theoretical Physics \\
Metrologichna str. 14-B\\
Kiev 03143, Ukraine

\end{center}

\vspace{1 cm}

\begin{abstract}
We describe properties of the nonstandard $q$-deformation
$U'_q({\rm so}_n)$ of the universal
enveloping algebra $U({\rm so}_n)$ of the Lie algebra ${\rm so}_n$
which does not coincide with the Drinfeld--Jimbo quantum algebra
$U_q({\rm so}_n)$. Irreducible representations of this algebras
for $q$ a root of unity $q^p=1$ are given. These representations act
on $p^N$-dimensional linear space (where $N$ is a number of positive
roots of the Lie algebra ${\rm so}_n$) and are given by
$r={\rm dim}\, {\rm so}_n$ complex parameters.
\end{abstract}

\vspace{1 cm}

We consider a $q$-deformation $U'_q({\rm so}_n)$ of the universal
enveloping algebra $U({\rm so}_n)$ of the Lie algebra ${\rm so}_n$
which does not coincide with the Drinfeld--Jimbo quantum algebra
$U_q({\rm so}_n)$. The algebra $U'_q({\rm so}_n)$ is constructed without
using the Cartan subalgebra and roots. This algebra has no root
elements.

The Drinfeld--Jimbo algebra $U_q({\rm so}_n)$ is obtained by deforming
Serre's relations for generating elements $E_1,\cdots , E_l$,
$F_1,\cdots , F_l$, $H_1,\cdots , H_l$ of $U_q({\rm so}_n)$. In order
to obtain $U'_q({\rm so}_n)$ we have to take determining relations for
the generating elements $I_{21}$, $I_{32},\cdots ,I_{n,n-1}$ of
$U_q({\rm so}_n)$ (they do not coincide with $E_j$, $F_j$, $H_j$) and to
deform these relations. The elements $I_{21}$, $I_{32},\cdots ,I_{n,n-1}$
belong to the basis $I_{ij}$, $i>j$, of the Lie algebra ${\rm so}_n$.
The matries $I_{ij}$, $i>j$, are defined as $I_{ij}=E_{ij}-E_{ji}$,
where $E_{ij}$ is the matrix with entries $(E_{ij})_{rs}=
\delta _{ir}\delta _{js}$. The universal enveloping algebra
$U({\rm so}_n)$ is generated by a part of the basis elements $I_{ij}$,
$i>j$, namely, by the elements $I_{21}$, $I_{32},\cdots ,I_{n,n-1}$.
These elements satisfy the relations
$$
I^2_{i,i-1}I_{i+1,i}-2I_{i,i-1}I_{i+1,i}I_{i,i-1} +
I_{i+1,i}I^2_{i,i-1} =-I_{i+1,i}, $$
$$
I_{i,i-1}I^2_{i+1,i}-2I_{i+1,i}I_{i,i-1}I_{i+1,i} +
I^2_{i+1,i}I_{i,i-1} =-I_{i,i-1}, $$
$$
I_{i,i-1}I_{j,j-1}- I_{j,j-1}I_{i,i-1}=0\ \ \ \ {\rm for}\ \ \ \
|i-j|>1.
$$
\noindent
{\bf Theorem 1.} {\it The universal enveloping algebra $U({\rm so}_n)$
is isomorphic to the complex associative algebra (with a unit element)
generated by the elements $I_{21}$, $I_{32},\cdots ,I_{n,n-1}$
satisfying the above relations.}
\medskip

We make the $q$-deformation of these relations by
$2\to [2]:=(q^2-q^{-2})/(q-q^{-1})=q+q^{-1}$.
As a result, we obtain the complex associative algebra
generated by elements $I_{21}$, $I_{32},\cdots ,I_{n,n-1}$ satisfying
the relations
$$
I^2_{i,i-1}I_{i+1,i}-(q+q^{-1})I_{i,i-1}I_{i+1,i}I_{i,i-1} +
I_{i+1,i}I^2_{i,i-1} =-I_{i+1,i}, $$
$$
I_{i,i-1}I^2_{i+1,i}-(q+q^{-1})I_{i+1,i}I_{i,i-1}I_{i+1,i} +
I^2_{i+1,i}I_{i,i-1} =-I_{i,i-1}, $$
$$
I_{i,i-1}I_{j,j-1}- I_{j,j-1}I_{i,i-1}=0\ \ \ \ {\rm for}\ \ \ \
|i-j|>1.
$$
This algebra was introduce by us in [1] and is denoted by
$U'_q({\rm so}_n)$.
There are the following motivations for studying this algebra and its
representations:
\medskip

1. The algebra of observables in 2+1 quantum gravity is isomorphic
to the algebra $U'_q({\rm so}_n)$ or to its quotient algebra [2, 3].
\smallskip

2. A quantum analogue of the Riemannian symmetric space $SU(n)/SO(n)$
is constructed by means of the algebra $U'_q({\rm so}_n)$ (see [4]).
\smallskip

3. A $q$-analogue of the theory of harmonic polynomials ($q$-harmonic
polynomials on quantum vector space ${\bf C}_q^n$) is constructed by
using the algebra $U'_q({\rm so}_n)$ (see [5]).
\medskip

In $U'_q({\rm so}_n)$ we can determine [6] elements analogous to the
matrices $I_{ij}$, $i>j$, of the Lie algebra ${\rm so}_n$. In order
to give them we use the
notation $I_{k,k-1}\equiv I^+_{k,k-1}\equiv I^-_{k,k-1}$.
Then for $k>l+1$ we define recursively
$$
I^{\pm}_{kl}:= [I_{l+1,l},I_{k,l+1}]_{q^{\pm 1}}\equiv
q^{\pm 1/2}I_{l+1,l}I_{k,l+1}-
q^{\mp 1/2}I_{k,l+1}I_{l+1,l}.
$$
(Note that similar sets of elements of $U'_q({\rm so}_n)$ are also
introduced in [5]).
The elements $I^+_{kl}$, $k>l$, satisfy the commutation relations
$$
[I^+_{ln},I^+_{kl}]_q=I^+_{kn},\ \
[I^+_{kl},I^+_{kn}]_q=I^+_{ln},\ \
[I^+_{kn},I^+_{ln}]_q=I^+_{kl} \ \ \
{\rm for}\ \ \  k>l>n,  $$
$$
[I^+_{kl},I^+_{nr}]=0\ \ \ \ {\rm for}\ \ \
k>l>n>r\ \ {\rm and}\ \ k>n>r>l,
$$
$$
[I^+_{kl},I^+_{nr}]_q=(q-q^{-1})
(I^+_{lr}I^+_{kn}-I^+_{kr}I^+_{nl}) \ \ \ {\rm for}\ \ \
k>n>l>r.
$$
For $I^-_{kl}$, $k>l$, the commutation relations are obtained by replacing
$I^+_{kl}$ by $I^-_{kl}$ and $q$ by $q^{-1}$.

Using the diamond lemma (see, for example, Chapter 4 in [7]),
N. Iorgov proved the Poincar\'e--Birkhoff--Witt theorem for the
algebra $U'_q({\rm so}_n)$ (the proof will be published):
\medskip

\noindent
{\bf Theorem 2.} {\it The elements ${I_{21}^+}^{m_{21}}
{I_{31}^+}^{m_{31}} {I_{32}^+}^{m_{32}} \cdots
{I_{n1}^+}^{m_{n1}}\cdots {I_{n,n-1}^+}^{m_{n,n-1}}$,
$m_{ij}=0,1,2$, $\cdots $,
form a basis of the algebra $U'_q({\rm so}_n)$.
This assertion is true if $I^+_{ij}$ are replaced by the
corresponding elements $I^-_{ij}$.}
\medskip

The algebra $U'_q({\rm so}_n)$ can be embedded into the Drinfeld--Jimbo
quantum algebra $U_q({\rm sl}_n)$. This quantum algebra is generated by
the elements $E_i$, $F_i$, $K_i^{\pm 1}=q^{\pm H_i}$,
$i=1,2,\cdots ,n-1$. Let us introduce the elements
$$
{\tilde I}_{j,j-1}=F_{j-1}-qq^{-H_{j-1}}E_{j-1},\ \ \ \  j=2,3,\cdots ,n.
$$

\noindent
{\bf Theorem 3.} {\it The homomorphism
$\varphi : U'_q({\rm so}_n)\to U_q({\rm sl}_n)$
uniquely determined by the relations
$\varphi (I_{i+1,i})={\tilde I}_{i+1,i}$, $i=1,\cdots ,n-1$, is an
isomorphism of $U'_q({\rm so}_n)$ to $U_q({\rm sl}_n)$.}
\medskip

That the mapping
$\varphi : U'_q({\rm so}_n)\to U_q({\rm sl}_n)$
is a homomorphism is proved in [4]. In [5] the authors state that it is
an isomorphism and say that it can be proved by means of the diamond lemma.
However, we could not restore their proof and found another proof.
It use the above Poincar\'e--Birkhoff--Witt theorem for the
algebra $U'_q({\rm so}_n)$. Namely, we use the explicit expressions for the
elements $I_{ij}\in U'_q({\rm so}_n)$ in terms of the elements of
the $L$-functionals of the quantum algebra $U_q({\rm sl}_n)$
from Subsec. 8.5.2 in [7]. Then we take any linear combination $a$
of the basis elements of $U'_q({\rm so}_n)$ from
Theorem 2 and express each element $I_{ij}$ in $a$ in terms of
elements of $U_q({\rm sl}_n)$. In this way we obtain the element $a$
expressed in term of elements of $U_q({\rm sl}_n)$. We consider it as a
linear combination of basis elements of $U_q({\rm sl}_n)$ given by
the Poincar\'e--Birkhoff--Witt theorem for $U_q({\rm sl}_n)$.
Next we show (by using the gradation of $U_q({\rm sl}_n)$ from Subsec.
6.1.5 in [7]) that at least one coefficient in this linear combination
is nonvanishing. This proves that
$\varphi$ is an isomorphism from $U'_q({\rm so}_n)$ to $U_q({\rm sl}_n)$.
\medskip

\noindent
{\bf Corollary.} {\it Let $q$ be not a root of unity. Then}

(a) {\it Finite dimensional irreducible representations
of $U'_q({\rm so}_n)$ separate elements of this algebra.}

(b) {\it The homomorphism $\psi : U'_q({\rm so}_3)\to {\hat U}_q({\rm
sl}_2)$
of the algebra $U'_q({\rm so}_3)$ to the extention ${\hat U}_q({\rm sl}_2)$
of the quantum algebra $U_q({\rm sl}_2)$ (given in [8]) uniquely
determined by the relations $\psi (I_{21})=\frac{\rm i}{q-q^{-1}}
(q^H-q^{-H})$ and $\psi (I_{32})=(E-F)(q^H+q^{-H})^{-1}$
is injective.}
\medskip

The assertion (a) follows from Theorem 3 and from the theorem on
separation of elements of the algebra $U_q({\rm sl}_n)$ by its
representations (see Subsec. 7.1.5 in [7]).
In order to prove (b), we consider representations of $U'_q({\rm so}_3)$
which are obtained by the composition $T\circ \psi$, where $T$ are
irreducible finite dimensional representations of ${\hat U}_q({\rm sl}_2)$.
It is shown in [8] that every irreducible finite dimensional
representation can be obtained from this composition. Now the assertion
(b) follows from (a).
\medskip

Now let us consider finite dimensional irreducible representations of
the algebra $U'_q({\rm so}_n)$. If $q$ is not a root of unity, there
exists two types of such representations: (a) representations of
the classical type (at $q\to 1$ they give the corresponding
finite dimensional irreducible representations of
the Lie algebra ${\rm so}_n$); (b) representations of the nonclassical
type (they do not admit the limit $q\to 1$ since in this point the
representation operators are singular). These representations are
described in [9].

We give irreducible representations of $U'_q({\rm so}_n)$ for
$q$ a root of unity. Using a description of central elements
of the algebra $U'_q({\rm so}_n)$ for $q$ a root of unity [10] it is proved
that every irreducible representation of $U'_q({\rm so}_n)$ in this
case is finite dimensional.

Let $q^k=1$ and $k$ is a smallest positive integer with this property.
We fix complex numbers $m_{1,n}, m_{2,n},..., m_{\left
\{ {n/2}\right \} ,n}$ (here $\{ {n/2}\}$ denotes integral part
of ${n/2}$) and $c_{ij}$, $h_{ij}$, $j=2,3,\cdots ,n-1$,
$i=1,2,\cdots , \{ {j/2}\}$ such that the differences
$h_{ij}-h_{sj}$ and the sums $h_{ij}+h_{sj}$ are not integers.
(We also suppose that $h_{p,2p+1}$ are not half-integral.)
The set of these numbers will be denoted by $\omega$.
Let $V$ be a complex vector space with a basis labelled by the
tableaux
$$
  \{\xi_{n} \}
\equiv \left\{ \matrix{ {\bf m}_{n} \cr {\bf m}_{n-1} \cr \dots
\cr {\bf m}_{2}  }
 \right\}
\equiv \{ {\bf m}_{n},\xi_{n-1}\}\equiv \{{\bf m}_{n} ,
{\bf m}_{n-1} ,\xi_{n-2}\} ,
                                                       \eqno(1)
$$
where the set of numbers ${\bf m}_{n}$
consists of $\{ {n/2}\}$ numbers $m_{1,n}, m_{2,n},\cdots , m_{\left
\{ {n/2}\right \} ,n}$ given above, and for each $s=2,3,\cdots ,n-1$,
${\bf m}_{s}$ is a set of numbers $m_{1,s}, \cdots , m_{\left
\{ {s/2}\right \} ,s}$ and each $m_{i,s}$ runs independently the
values $h_{i,s}, h_{i,s}+1,\cdots , h_{i,s}+k-1$.
Thus, ${\rm dim}\, V$ coincides with $k^N$, where $N$ is the number of
positive roots of ${\rm so}_n$.
It is convenient to use for the numbers $m_{i,s}$ the so-called
$l$-coordinates
$$
l_{j,2p+1}=m_{j,2p+1}+p-j+1,  \qquad
                      l_{j,2p}=m_{j,2p}+p-j .        \eqno (2)
$$

To the set of numbers $\omega$ there corresponds the
irreducible finite dimensional representation $T_\omega$ of the algebra
$U'_q({\rm so}_n)$.
The operator $T_{\omega}(I_{2p+1,2p})$ of the representation
$T_{\omega}$ acts upon the basis elements, labelled by (1), by the formula
$$
T_{\omega}(I_{2p+1,2p})
| \xi_n\rangle =
\sum^p_{j=1} c_{j,2p} \frac{ A^j_{2p}(\xi_n)}
{q^{l_{j,2p}}+q^{-l_{j,2p}} } \,
            \vert (\xi_n)^{+j}_{2p}\rangle -
\sum^p_{j=1} c^{-1}_{j,2p} \frac{A^j_{2p}((\xi_n)^{-j}_{2p})}
{q^{l_{j,2p}}+q^{-l_{j,2p}}} \,
|(\xi_n)^{-j}_{2p}\rangle                         \eqno(3)
$$
and the operator $T_{\omega}(I_{2p,2p-1})$ of the representation
$T_{\omega}$ acts as
$$
T_{\omega}(I_{2p,2p-1})\vert \xi_n\rangle=
\sum^{p-1}_{j=1}c_{j,2p-1}  \frac{B^j_{2p-1}(\xi_n)}
{[2 l_{j,2p-1}-1][l_{j,2p-1}]} \,
\vert (\xi_n)^{+j}_{2p-1} \rangle
$$
$$
-\sum^{p-1}_{j=1}c^{-1}_{j,j,2p-1} \frac {B^j_{2p-1}((\xi_n)^{-j}_{2p-1})}
{[2 l_{j,2p-1}-1][l_{j,2p-1}-1]} \,
\vert (\xi_n)^{-j}_{2p-1}\rangle
+ {\rm i}\, C_{2p-1}(\xi_n) \,
\vert \xi_n \rangle ,                                 \eqno(4)
$$
where numbers in square brackets mean $q$-numbers:
$[a]:= (q^a-q^{-a})/(q-q^{-1})$.

In these formulas, $(\xi_n)^{\pm j}_{s}$ means the tableau (1)
in which $j$-th component $m_{j,s}$ in ${\bf m}_s$ is replaced
by $m_{j,s}\pm 1$. If $m_{j,s}+1=h_{j,s}+k$ (resp.
$m_{j,s}-1=h_{j,s}-1$), then we set $m_{j,s}+1=h_{j,s}$ (resp.
$m_{j,s}-1=h_{j,s}+k-1$). The coefficients
$A^j_{2p},  $ $B^j_{2p-1},$ $C_{2p-1}$ in (3) and (4) are given
by the expressions
$$
A^j_{2p}(\xi_n) {=}
\left( \frac{\prod_{i=1}^p [l_{i,2p+1}+l_{j,2p}] [l_{i,2p+1}{-}l_{j,2p}{-}1]
\prod_{i=1}^{p-1} [l_{i,2p-1}+l_{j,2p}] [l_{i,2p-1}{-}l_{j,2p}{-}1]}
{\prod_{i\ne j}^p [l_{i,2p}+l_{j,2p}][l_{i,2p}-l_{j,2p}]
[l_{i,2p}+l_{j,2p}+1][l_{i,2p}-l_{j,2p}-1]} \right)^{1\over 2} ,
$$
$$
B^j_{2p-1}(\xi_n){=}  \! \!
\left( \frac{\prod_{i=1}^p
[l_{i,2p}+l_{j,2p-1}] [l_{i,2p}-l_{j,2p-1}] \prod_{i=1}^{p-1}
[l_{i,2p-2}+l_{j,2p-1}] [l_{i,2p-2}-l_{j,2p-1}]}
{\prod_{i\ne j}^{p-1}
[l_{i,2p-1}{+}l_{j,2p-1}][l_{i,2p-1}{-}l_{j,2p{-}1}]
[l_{i,2p-1}{+}l_{j,2p-1}{-}1][l_{i,2p-1}{-}l_{j,2p-1}{-}1]}
 \right) ^{1\over 2}
$$
$$
C_{2p-1}(\xi_n) =\frac{ \prod_{s=1}^p [ l_{s,2p} ]
\prod_{s=1}^{p-1} [ l_{s,2p-2} ]}
{\prod_{s=1}^{p-1} [l_{s,2p-1}] [l_{s,2p-1} - 1] } .
$$
It is seen from (2) that $C_{2p-1}(\xi_n)=0$
if $m_{p,2p}\equiv l_{p,2p}=0$.

The representations $T_\omega$ are given by
$r={\rm dim}\, {\rm so}_n$ complex numbers.
These representations do not exhaust all irreducible
representations of $U'_q({\rm so}_n)$ for $q$ a root of unity.
The other irreducible representations as well as a proof of
the fact that the above formulas determine
representations of $U'_q({\rm so}_n)$ will be published.

\section*{Acknowledgment}
The research of this paper was made possible in part by Award
UP1-309 of Civilian Research and Development Foundation and by
Award 1.4/206 of Ukrainian DFFD.


\begin{thebibliography}{**}
\bibitem{}A. M. Gavrilik and A. U. Klimyk, {\it Lett. Math. Phys.} {\bf 21}
(1991), 215.
\bibitem{}J. Nelson, T. Regge, and F. Zertuche, {\it Nucl. Phys.} {\bf B339}
(1990), 227.
\bibitem{}J. Nelson and T. Regge, {\it Commun. Math. Phys.} {\bf 155}
(1993), 561.
\bibitem{}M. Noumi, {\it Adv. Math.} {\bf 123} (1993), 16.
\bibitem{}M. Noumi, T. Umeda, and M. Wakayama, {\it Compos. Math.}
{\bf 104} (1996), 227.
\bibitem{}A. M. Gavrilik and N. Z. Iorgov, {\it Ukr. J. Phys.} {\bf 43}
(1998), 456.
\bibitem{} A. Klimyk and K. Schm\"udgen: {\it Quantum Groups and
Their Representations}, Springer, Berlin, 1997.
\bibitem{} M. Havl\'\i \v cek, A. Klimyk, and S. Po\v sta,
{\it J. Math. Phys.} {\bf 40} (1999), 2135.
\bibitem{}N. Z. Iorgov and A. U. Klimyk, {\it Czech. J. Phys.} {\bf 49}
(1999), to be published.
\bibitem{} M. Havl\'\i \v cek, A. Klimyk, and S. Po\v sta,
{\it Czech. J. Phys.} {\bf 49} (1999), to be published.






\end{thebibliography}
\end{document}